\documentclass[draft,a4paper,11pt]{article}
\newtheorem{theorem}{Theorem}[section]
\newtheorem{lemma}{Lemma}[section]

\newtheorem{proposition}{Proposition}[section]

\newenvironment{proof}{{\bf Proof.}}{\par\hspace{25em}\rule{1ex}{1ex}\par}
\newcommand{\Lin}{\mathop{\rm Lin\,}}
\newcommand{\Cal}[1]{{\cal #1}}
\newcommand{\pfend}{\par\hspace{25em}\rule{1ex}{1ex}\par}
\title{Powers of the space forms  curvature operator and geodesics of
       the tangent bundle.\footnote{Ukr. Math. Journal  2004, 56/9,  1231-1243.}}
\author{Saharova Y., Yampolsky A.}
\date{}
\begin{document}
\maketitle
\begin{abstract}
    It is well-known that if $\Gamma$ is a geodesic line of the tangent (sphere)
    bundle with Sasaki metric of a locally symmetric Riemannian manifold then the
    projected curve $\gamma=\pi\circ\Gamma$ has all its geodesic curvatures constant. In
    this paper we consider the case of tangent (sphere) bundle
    over the real, complex and quaternionic space form and give a unified proof of the
    following property: all geodesic curvatures of projected curve are zero starting
    from $k_3,\,k_6$ and $k_{10}$ for the real, complex and quaternionic space formes
    respectively.
    \\[2ex]
    {\it AMS subject class:} Primary 54C40,14E20; Secondary 46E25, 20C20
\end{abstract}
\section*{Introduction}

Sato K. \cite{St}  and Sasaki S. \cite{Sk} proved that the
projection to the base space of any non-vertical geodesic line on
the tangent or the tangent sphere bundle of a real space form
$M^n(c)$ is a curve of constant curvatures $k_1$ and $k_2$ and
zero curvatures $k_3,\dots, k_{n-1}$. Nagy P. \cite{Ng}
essentially generalized this result. He considered the case of
general locally symmetric base manifold and have proved that the
geodesic curvatures of projection of any (non-vertical) geodesic
line on the tangent sphere bundle are all constant. Nevertheless
it was still interesting to find a clearer description of
projections of geodesics for the case of classical rank one
symmetric spaces. The second author made a first step in this
defection and proved that the projection to the base space of any
non-vertical geodesic line on the tangent or tangent sphere bundle
of a complex space form $CP^n$ is a curve of constant curvatures
$k_1, \dots, k_5$  and zero curvatures $k_6,\dots, k_{n-1}$

In is this paper we make a contribution in more clear
understanding of geometry of projected geodesics in the case of
tangent (sphere) bundle of  almost all classical locally symmetric
spaces, namely \textit{spheres, complex and quaternionic
projective spaces and their non-compact dual} from a unified
viewpoint using the \textit{recurrent properties} of powers of the
curvature operator of these spaces. This approach allows to give
also a unified proof of the results from \cite{Sk}, \cite{St} and
\cite{Ym}

We also use an easy to prove result \cite{Azo}, stating that the
geodesics of tangent or tangent sphere bundle with Sasaki metric
have the same projections to the base manifold.

    \textbf{Remark on notations.} Throughout the paper $\big<\cdot,\cdot\big>$ and $|\cdot|$
    mean the scalar product and the norm  of vectors with respect to the corresponding
    metrics.

\section{Summary of main results.}

    Let $(M^n(c),g)$ be a Riemannian manifold of constant curvature $c$,
    \\ $(M^{2n}(c);J;g)$  a Riemannian manifold with complex structure $J$ of
    constant holomorphic curvature $c$ and $(M^{4n}(c);J_1,J_2,J_3;g)$ a Riemannian manifold with
    quaternionic structure $(J_1,J_2,J_3)$ of constant quaternionic curvature $c$. For
    the sake of brevity, denote by $\Cal{M}(c)$ one of these space forms with corresponding
    standard metrics and will refer to
    $\Cal{M}(c)$ just to a space  form of constant curvature $c$.
    The main result is the following statement.

    \begin{theorem} \label{Main}
    Let $\Cal{M}(c)$ be a space form of constant curvature $c\ne0$. Let $\Gamma$ be
    non-vertical geodesic line on the tangent or tangent sphere bundle over $\mathcal{M}(c)$. Let
    $\gamma=\pi\circ\Gamma$ be the projection of $\Gamma$ to $\Cal{M}(c)$. Then
    the geodesic curvatures $k_1,k_2,\dots$ of $\gamma$ are all constant and
    \begin{itemize}\itemsep=-1ex
    \item[\rm (a)] $k_3=\dots =k_{n-1}=0$ for the real space form;
    \item[\rm (b)] $k_6=\dots=k_{2n-1}=0$ for the complex space form;
    \item[\rm (c)] $k_{10}=\dots=k_{4n-1}=0$ for the quaternionic space form.
    \end{itemize}
    \end{theorem}

    As the referee remarked, the result of the Theorem \ref{Main} can be
    expressed in more clear geometrical terms, namely
    \textit{ the projected curve $\gamma=\pi\circ\Gamma$ lies in a
    totally geodesic $S^3$ or $H^3$, in a totally geodesic $CP^3$ or
    $CH^3$ and in a totally geodesic $QP^3$ or $QH^3$ for the real,
    complex and quaternionic space form respectively.} These
    assertions can be derived from (\ref{tgS3}), (\ref{tgCP3}) and
    (\ref{tgCH3}).

    Proof of the Theorem \ref{Main} is based on the recurrent property of powers of curvature operator of spaces under
    consideration. Let $R_{XY}$ be the curvature operator of $\Cal{M}(c)$.
    Define a power of curvature operator $R^p_{XY}$ recurrently in the following way:
$$
R^{p}_{XY}Z=R^{p-1}_{XY}(R_{XY}Z)\ \ p>1.
$$
    The basic tool for our considerations are a chain of lemmas.
\begin{lemma}\label{Const}
    Let $R_{XY}$ be the curvature operator of the real space form \\ $(M^n(c),g)$. Then for any
    $X$ and $Y$
    $$
    R^p_{XY}=
    \left\{
    \begin{array}{l}
    (-b^2c^2)^{s-1} R_{XY} \mbox{  for p=2s-1}\\[1ex]
     (-b^2c^2)^{s-1} R_{XY}^2 \mbox{  for p=2s},
    \end{array}
    \right.
    s\geq1
    $$
    where $b=|X\wedge Y|$ is a norm of bivector $X\wedge Y$.
\end{lemma}

\begin{lemma}\label{Comp}
    Let $R_{XY}$ be the curvature operator of the non-flat complex space form $(M^n(c);J;g)$.
    Denote by $b=|X\wedge Y|$ the norm of a bivector $X\wedge Y$ and $m=\big<X,JY\big>$.
    Then for any $X$ and $Y$
$$
    R^p_{XY}=
    \left\{
    \begin{array}{l}
    \Lin(JR^2_{XY},\,R_{XY},\ J\,) \mbox{  for p=2s-1}\\[1ex]
    \Lin(R^2_{XY},\,JR_{XY},\,E\,) \mbox{  for p=2s},
    \end{array}
    \right.
    s\geq2
$$
    where $E$ is the identity operator and $\Lin$ means a linear combination of corresponding
    operators with coefficients being polynomials in \,$\frac1c\,, b\,,m$.
\end{lemma}

\begin{lemma}\label{Quat}
    Let $R_{XY}$ be the curvature operator of the non-flat quaternionic space form
     $(M^n(c);J_1,J_2,J_3;g)$.
    Denote by $b=|X\wedge Y|$ the norm of a bivector $X\wedge Y$. Set  $m_1=\big<X,J_1Y\big>$,
    $m_2=\big<X,J_2Y\big>$, $m_3=\big<X,J_3Y\big>$, $m^2=m_1^2+m_2^2+m_3^2$,
    ${\cal J}=m_1J_1+m_2J_2+m_3J_3$.
    Then for any $X$ and $Y$
$$
    R^p_{XY}=
    \left\{
    \begin{array}{l}
    \Lin({\cal J}R^4_{XY},\,R^3_{XY},\,{\cal J}R_{XY}^2,\,R_{XY}, {\cal J}\,) \mbox{  for
    p=2s-1}\\[1ex]
    \Lin(R^4_{XY},\,{\cal J}R^3_{XY},\,R_{XY}^2,\, {\cal J}R_{XY},E\,) \mbox{  for p=2s},
    \end{array}
    \right.
    s\geq3
$$
    where $E$ is the identity operator and $\Lin$ means a linear combination of corresponding
    operators with coefficients  being polynomials in \,$\frac1c\,, b\,,m$.
\end{lemma}

\section{Necessary facts and proof of the main result.}

    Let $(M^n,g)$ be a Riemannian manifold and $TM^n$ be its tangent bundle. Denote by
    $(u^1,\dots,u^n)$ a local coordinate system on $M^n$. Then in each tangent space of $M^n$
    the natural coordinate frame $\big\{ \partial/\partial u^1,\dots,\partial/\partial u^n\big\}$
    form a local basis. Let $\xi$ be any tangent vector over the given local chart.
    Then $\xi$ can be decomposed as
$$
    \xi=\xi^1 \partial/\partial u^1+\dots+\xi^n\partial/\partial u^n.
$$
    The parameters $(u^1,\dots,u^n;\xi^1,\dots,\xi^n)$ form the so called
    \textit{natural induced coordinate system} in $TM^n$. The \textit{Sasaki metric} line
    element $d\sigma ^2$ with respect to this coordinate system is
\begin{equation}\label{mtc}
    d\sigma^2=ds^2+|D\xi|^2,
\end{equation}
    where $ds^2$ is a line element of $M^n$, $D\xi$ is the covariant differential of
    $\xi$ with respect to Levi-Civita  connection on $M^n$ and $|\cdot|$ means the
    norm with respect to Riemannian metric on $M^n$.

    The \textit{tangent sphere bundle } $T_1M^n$ can be considered as a hypersurface in
    the tangent bundle defined by the condition $|\,\xi\,|=1$. We will consider
    $T_1M^n$ as a submanifold in $TM^n$ with the induced metric.

    With respect to the natural coordinate system, each curve $\Gamma$ on $TM^n$ can be
    represented as $\Gamma(\sigma)=\Big\{u^1(\sigma)\dots,u^n(\sigma);
    \xi^1(\sigma),\dots,\xi^n(\sigma)\Big\}$
    with respect to the arc-length parameter $\sigma$ and can be interpreted as the
    vector field $\xi(\sigma)=\xi^1(\sigma)\partial/\partial u^1+
    \dots+\xi^n(\sigma)\partial/\partial u^n$ along the \textit{projected} curve
    $\gamma=\pi\circ\Gamma=(u^1(\sigma),\dots,u^n(\sigma))$. If $\xi$ is a \textit{ unit}
    vector field then $\Gamma$ lies in $T_1M^n$ and represents an arbitrary curve in
    $T_1M^n$.

    Denote by $(')$ the covariant derivative along $\gamma$ with respect to parameter
    $\sigma$. Then $\Gamma$ is a geodesic line on $TM^n$ or $T_1M^n$ if $\gamma$ and
    $\xi$ satisfy respectively the system of equations
    \begin{equation}\label{EqGeo}
    \begin{array}{c}
       \mbox{$TM^n:$  } \left\{\begin{array}{l}
                \gamma''=R_{\,\xi'\,\xi\,} \gamma',\\[1ex]
                \xi''=0;
                \end{array}
        \right.
        \qquad
        \mbox{$T_1M^n:$  }\left\{\begin{array}{l}
                \gamma''=R_{\,\xi'\,\xi\,} \gamma',\\[1ex]
                \xi''=-\rho^2\xi,
                \end{array}
        \right.
    \end{array}
    \end{equation}
    where $\rho^2=|\,\xi'\,|^2$ and $R_{\,\xi'\,\xi\,}$ is the \textit{curvature operator}
    of $M^n$ based on bivector $\xi'\wedge\xi$.

    From (\ref{EqGeo}) it follows that $\rho=const$  in both cases.  Denote by $s$ the
    arc-length parameter on $\gamma$. Then from (\ref{mtc}) it follows that
    \begin{equation}\label{dsdsig}
    \frac{ds}{d\sigma}=\sqrt{1-\rho^2},
    \end{equation}
    so that $0\leq \rho\leq 1$. According to the latter inequality, the set of geodesics of $TM^n$
    and $T_1M^n$ can be splitted naturally into 3 classes, namely,
    \begin{itemize}\itemsep=-1ex
    \item  \textit{horizontal} geodesics ($\rho=0$) generated by parallel (unit) vector fields along the
    geodesics on the base manifold;
    \item \textit{vertical} geodesics ($\rho=1$) represented by geodesics on a fixed
    fiber;
    \item \textit{umbilical} geodesics corresponding to $0<\rho<1$.
    \end{itemize}

    In what follows we will consider the properties of projections of umbilical
    geodesics.
    \begin{lemma}(cf. \cite{Ng})\label{deriv}
    Let $(M^n,g)$ be a locally symmetric Riemannian manifold and $R_{XY}$ its curvature
    operator. Let $\gamma=\pi\circ\Gamma$ be a projection of geodesic line on $TM^n$ or
    $T_1M^n$ to the base space. Then for the derivatives of $\gamma$ of order $p$ we have
    $$
    \gamma^{(p)}=R_{\,\xi'\,\xi}^{p-1}\gamma'=R_{\,\xi'\,\xi}\gamma^{(p-1)}
    $$
    and as a consequence all the geodesic curvatures of $\gamma$ are constant.
   \end{lemma}
   \begin{proof}
   The equalities follow from parallelism of curvature tensor of $M^n$ and the
   equations (\ref{EqGeo}). Moreover, from the evident identity
   $$
   \big<\gamma^{(p)},\gamma^{(p-1)}\big>\equiv0
   $$
   for all $p>1$, we conclude that $|\gamma^{(p)}|=const$ for all $p>1$ and
   therefore, by induction, all the geodesic curvatures of $\gamma$ are constant.
   \end{proof}

   \textbf{Proof of Theorem \ref{Main}.}
   \textbf{Case (a).}  Denote by $e_1,\dots,e_{n-1}$ the Frenet frame of $\gamma$. Using the Frenet
   formulas for the curve with constant geodesic curvatures and keeping in mind (\ref{dsdsig}),
   it is easy to see that
   \begin{equation}\label{frenet}
   \begin{array}{l}
   \gamma^{(2s-1)}=(1-\rho^2)^{s-1/2}k_1k_2\dots k_{2s-2}\,e_{2s-1} + \Lin
   \big\{e_1,e_3,\dots,e_{2s-3}\big\},\\[1ex]
   \gamma^{(2s)}\quad=(1-\rho^2)^{s}k_1k_2\dots k_{2s-1}\,e_{2s} + \Lin
   \big\{e_2,e_4,\dots,e_{2s-2}\big\}
   \end{array}
   \end{equation}
   for all $s\geq1$ ( with formal setting $k_0\equiv1$). Setting $s=1,2$ in even derivatives,
   we see that
   \begin{equation}\label{derivR}
   \begin{array}{l}
   \gamma^{(2)}=(1-\rho^2)k_1\,e_2\\[1ex]
   \gamma^{(4)}=(1-\rho^2)^2k_1k_2k_3\,e_4+\Lin (e_2).
   \end{array}
   \end{equation}
   On the other hand, applying Lemma \ref{deriv}, Lemma \ref{Const} and Lemma
   \ref{deriv} again, we get
   \begin{equation}\label{tgS3}
   \gamma^{(4)}=R_{\,\xi'\,\xi}^3\gamma'=-b^2c^2\,R_{\,\xi'\,\xi}\gamma'=-b^2c^2\gamma^{(2)}.
   \end{equation}
   Substitution from (\ref{derivR}) gives
   $$
   (1-\rho^2)^2k_1k_2k_3\,e_4+\Lin (e_2)=0
   $$
   and therefore $k_3=0$, which completes the proof.

   Remark, that $b^2$ is constant along $\gamma$ since
   $$
   \big(b^2\big)'=\Big(|\xi'\wedge\xi|^2\Big)=
   \Big(\rho^2\,|\,\xi\,|^2-\big<\xi',\xi\big>^2\Big)'=
   2\rho^2\big<\xi',\xi\big>-2\big<\xi',\xi\big>\rho^2\equiv0.
   $$

   \textbf{Case (b).} Denote by $e_1,\dots,e_{2n-1}$ the Frenet frame of $\gamma$. Similar to the case
   (a) considerations,  the Frenet formulas give
   \begin{equation}\label{frenetC}
   \begin{array}{l}
   \gamma^{(2s-1)}=(1-\rho^2)^{s-1/2}k_1k_2\dots k_{2s-2}\,e_{2s-1} + \Lin
   \big\{e_1,e_3,\dots,e_{2s-3}\big\},\\[1ex]
   \gamma^{(2s)}\quad=(1-\rho^2)^{s}k_1k_2\dots k_{2s-1}\,e_{2s} + \Lin
   \big\{e_2,e_4,\dots,e_{2s-2}\big\}
   \end{array}
   \end{equation}
   for all $s\geq1$. Setting $s=1,2,3,4$ in odd derivatives, we get
   \begin{equation}\label{derivC}
   \begin{array}{l}
   \gamma'\quad=(1-\rho^2)^{1/2}e_1,\\[1ex]
   \gamma^{(3)}=(1-\rho^2)^{3/2}k_1k_2\,e_3+\Lin (e_1),\\[1ex]
   \gamma^{(5)}=(1-\rho^2)^{5/2}k_1\dots k_4\,e_5+\Lin (e_1,e_3),\\[1ex]
   \gamma^{(7)}=(1-\rho^2)^{7/2}k_1\dots k_6\,e_7+\Lin (e_1,e_3,e_5).
   \end{array}
   \end{equation}
   On the other hand, applying Lemma \ref{deriv}, Lemma \ref{Comp} and Lemma
   \ref{deriv} again, we get
   \begin{equation}\label{derTensC}
    \left\{
   \begin{array}{l}
   \gamma^{(5)}=R^4_{\,\xi'\,\xi}\gamma'=
   \Lin (R^2_{\,\xi'\,\xi},J R_{\,\xi'\,\xi}, E)\gamma'=
   \Lin (\gamma^{(3)},J\gamma^{(2)}, \gamma'),\\[1ex]
   \gamma^{(7)}=R^6_{\,\xi'\,\xi}\gamma'=
   \Lin (R^2_{\,\xi'\,\xi},J R_{\,\xi'\,\xi}, E)\gamma'=
   \Lin (\gamma^{(3)},J\gamma^{(2)}, \gamma'),\\[1ex]
   \end{array}
   \right.
   \end{equation}
   Excluding $J\gamma^{(2)}$ from (\ref{derTensC}), we come to the equation
   \begin{equation}\label{tgCP3}
    \gamma^{(7)}=\Lin (\gamma^{(5)},\gamma^{(3)},\gamma').
   \end{equation}
   Substitution from (\ref{derivC}) imply
    $$
    (1-\rho^2)^{7/2}k_1\dots k_6\,e_7+\Lin (e_1,e_3,e_5)=0
    $$
    and we conclude that $k_6=0$ which completes the proof.

   Remark, that the coefficients of all linear combinations are constants.
   Indeed, by Lemma \ref{Comp} the coefficients are polynomials in
   $1/c, b=|\xi'\wedge\xi|$ and $m=\big<\xi',J\xi\big>$.
   The value $b$ is constant along $\gamma$ by the same reasons as in case (a).
   The value $m$ is constant along $\gamma$ since
   $$
   m'=\big<\xi',J\xi\big>'=\big<\xi'',J\xi\big>+\big<\xi',J\xi'\big>\equiv0.
   $$

   \textbf{Case (c).} Denote by $e_1,\dots,e_{4n-1}$ the Frenet frame of $\gamma$. As above,
    the Frenet formulas give
   \begin{equation}\label{frenetH}
   \begin{array}{l}
   \gamma^{(2s-1)}=(1-\rho^2)^{s-1/2}k_1k_2\dots k_{2s-2}\,e_{2s-1} + \Lin
   \big\{e_1,e_3,\dots,e_{2s-3}\big\},\\[1ex]
   \gamma^{(2s)}\quad=(1-\rho^2)^{s}k_1k_2\dots k_{2s-1}\,e_{2s} + \Lin
   \big\{e_2,e_4,\dots,e_{2s-2}\big\}
   \end{array}
   \end{equation}
   for all $s\geq1$. Setting $s=1,2,3,4,5,6$ in odd derivatives, we get
   \begin{equation}\label{derivH}
   \begin{array}{l}
   \gamma'\quad=(1-\rho^2)^{1/2}e_1,\\[1ex]
   \gamma^{(3)}=(1-\rho^2)^{3/2}k_1k_2\,e_3+\Lin (e_1),\\[1ex]
   \gamma^{(5)}=(1-\rho^2)^{5/2}k_1\dots k_4\,e_5+\Lin (e_1,e_3),\\[1ex]
   \gamma^{(7)}=(1-\rho^2)^{7/2}k_1\dots k_6\,e_7+\Lin (e_1,e_3,e_5),\\[1ex]
   \gamma^{(9)}=(1-\rho^2)^{9/2}k_1\dots k_8\,e_9+\Lin (e_1,e_3,e_5,e_7),\\[1ex]
   \gamma^{(11)}=(1-\rho^2)^{11/2}k_1\dots k_{10}\,e_{11}+\Lin (e_1,e_3,e_5,e_7,e_9).
   \end{array}
   \end{equation}
   Applying again Lemma \ref{deriv}, Lemma \ref{Quat} and then Lemma
   \ref{deriv}, we get
   \begin{equation}\label{derTensH}
    \left\{
   \begin{array}{ll}
   \gamma^{(7)}=R^6_{\,\xi'\,\xi}\gamma'=&
   \Lin (R^4_{\,\xi'\,\xi},\Cal{J} R^3_{\,\xi'\,\xi},R^2_{\,\xi'\,\xi},\Cal{J}R_{\,\xi'\,\xi},
   E)\gamma'=\\[1ex]
   &\Lin (\gamma^{(5)},\Cal{J}\gamma^{(4)},\gamma^{(3)},\Cal{J}\gamma^{(2)}, \gamma'),\\[2ex]
   \gamma^{(9)}=R^8_{\,\xi'\,\xi}\gamma'=&
   \Lin (R^4_{\,\xi'\,\xi},\Cal{J} R^3_{\,\xi'\,\xi},R^2_{\,\xi'\,\xi},\Cal{J}R_{\,\xi'\,\xi},
   E)\gamma'=\\[1ex]
   &\Lin (\gamma^{(5)},\Cal{J}\gamma^{(4)},\gamma^{(3)},\Cal{J}\gamma^{(2)}, \gamma'),\\[2ex]
   \gamma^{(11)}=R^{10}_{\,\xi'\,\xi}\gamma'=&
   \Lin (R^4_{\,\xi'\,\xi},\Cal{J} R^3_{\,\xi'\,\xi},R^2_{\,\xi'\,\xi},\Cal{J}R_{\,\xi'\,\xi},
   E)\gamma'=\\[1ex]
   &\Lin (\gamma^{(5)},\Cal{J}\gamma^{(4)},\gamma^{(3)},\Cal{J}\gamma^{(2)}, \gamma').
   \end{array}
   \right.
   \end{equation}
   Excluding $\Cal{J}\gamma^{(2)}$  and $\Cal{J}\gamma^{(4)}$ from (\ref{derTensH}),
   we come to the equation
   \begin{equation}\label{tgCH3}
    \gamma^{(11)}=\Lin (\gamma^{(9)},\gamma^{(7)},\gamma^{(5)},\gamma^{(3)},\gamma').
   \end{equation}
   Substitution from (\ref{derivH}) imply
    $$
    (1-\rho^2)^{11/2}k_1\dots k_{10}\,e_{11}+\Lin (e_1,e_3,e_5,e_7,e_9)=0
    $$
    and we conclude that $k_{10}=0$ which completes the proof.

   Remark, that the coefficients of all linear combinations are constants.
   Indeed, by Lemma \ref{Quat} the coefficients are polynomials in
   $1/c, b=|\xi'\wedge\xi|$ and $m=\sqrt{m_1^2+m_2^2+m_3^2}$.
   The value $b$ is constant along $\gamma$ by the same reasons as in case (a).
   The values $m_1,m_2,m_3$ are all constant along $\gamma$ since
   $$
   m_i'=\big<\xi',J_i\xi\big>'=\big<\xi'',J_i\xi\big>+\big<\xi',J_i\xi'\big>\equiv0
   $$
   for $i=1,2,3$.

\section{ Proofs of basic Lemmas}
   {\bf Proof of Lemma \ref{Const}}.The curvature operator $R_{XY}$ of the real space form
   $(M^n(c),g)$ has the following expression
   $$
   R_{XY}Z=c\,\Big[\big<Y,Z\big>X-\big<X,Z\big>Y\Big].
   $$
   Then

   \noindent
   $
   R^2_{XY}Z=c\,\Big[\big<Y,R_{XY}Z\big>X-\big<X,R_{XY}Z\big>Y\Big]=
         c^2\Big[\Big<Y,\big<Y,Z\big>X-\big<X,Z\big>Y\Big>X-
         \Big<X,\big<Y,Z\big>X-\big<X,Z\big>Y\Big>Y\Big]=
         c^2\Big[\big<Y,Z\big>\big<X,Y\big>X-\big<X,Z\big>|Y|^2X-\\
         \big<Y,Z\big>|X|^2Y+
         \big<X,Z\big>\big<X,Y\big>Y\Big]=
          c^2\Big[\big<Y,Z\big>\Big(\big<X,Y\big>X-|X|^2Y\Big)+\\
          \big<X,Z\big>
         \Big(\big<X,Y\big>Y-|Y|^2X\Big)\Big]=
         c\,\Big[\big<Y,Z\big>R_{XY}X+\big<X,Z\big>R_{YX}Y\Big].
   $
\\
   Therefore,
\\
    \noindent
   $
   R^3_{XY}Z=c\Big[\big<Y,R_{XY}Z\big>R_{XY}X+\big<X,R_{XY}Z\big>R_{YX}Y\Big]=
c^3\Big[\Big(\big<Y,Z\big>\big<X,Y\big>-\big<X,Z\big>|Y|^2\Big)\Big(\big<X,Y\big>X-                                                                   |X|^2Y\Big)+
\Big(\big<Y,Z\big>|X|^2-\big<X,Z\big>\big<X,Y\big>\Big)                                            \Big(\big<X,Y\big>Y-|Y|^2X\Big)\Big]=
c^3\Big[-\big<Y,Z\big>X\Big(|X|^2|Y|^2-\big<X,Y\big>^2\Big)+
\big<X,Z\big>Y\Big(|X|^2|Y|^2-\\
 \big<X,Y\big>^2\Big)\Big]=-c^2b^2R_{XY}Z,
   $

\noindent
   where, evidently, $b^2=|X|^2|Y|^2-\big<X,Y\big>^2$ is the square norm of $X\wedge Y$.

   Now we can find the other powers for $R_{XY}$ inductively.
   \pfend

   \textbf{Proof of Lemma \ref{Comp}}

    The curvature operator $R_{XY}$ of the complex space form
   $(M^{2n}(c);J;g)$ has the following expression
    $$
   R_{XY}Z=\frac{c}{4}\Big[\big<Y,Z\big>X-\big<X,Z\big>Y+\big<JY,Z\big>JX-\big<JX,Z\big>JY+
   2\,\big<X,JY\big>JZ\Big].
   $$
   Introduce the unit sphere type operator $S$ acting as
   $$
   S(Z)\stackrel{def}{=} S_{XY}Z=\big<Y,Z\big>X-\big<X,Z\big>Y,
   $$
   and the operator $\hat S(Z)$ acting as
   $$
   \hat S(Z)\stackrel{def}{=} S_{JX\,JY}Z=\big<JY,Z\big>JX-\big<JX,Z\big>JY.
   $$
   Finally, if we denote $m=\big<X,JY\big>$, then the curvature operator under
   consideration takes the form
   \begin{equation}\label{curvop}
   R_{XY}Z=\frac{c}{4}\Big[ S+\hat S +2m\,J\Big]Z.
   \end{equation}
   Since $|X\wedge Y|=|(JX)\wedge (JY)|$, the operators $S$ and $\hat S$ satisfy
   $$
   S^3=-b^2S, \quad \hat S^3=-b^2\hat S,
   $$
   where $b^2=|X\wedge Y|^2$.

   In what follows, we need a "table of products" for the operators $S$ and $\hat S$.
   Namely,

   \begin{equation}\label{tableC}
   \begin{array}{|c|c|c|c|}
   \hline
   &S&\hat S&J\\[1ex]
   \hline S&S^2&mJ\hat S&J\hat S\\[1ex]
   \hline \hat S&mJS&\hat S^2&JS\\[1ex]
   \hline J&JS&J\hat S&-E\\[1ex]
   \hline
   \end{array}
   \end{equation}

   Indeed,
   $$
\begin{array}{l}
   \begin{array}{l}
   (S\hat S)(Z)=S_{XY}\hat S_{JX\,JY}Z=S_{XY}\Big[\big<JY,Z\big>JX
   -\big<JX,Z\big>JY\Big]=\\
   \Big<Y,\big<JY,Z\big>JX-\big<JX,Z\big>JY\Big>X-
    \Big<X,\big<JY,Z\big>JX-\big<JX,Z\big>JY\Big>Y=\\
    \big<Y,JX\big>\big<JY,Z\big>X+\big<JX,Z\big>\big<X,JY\big>Y=\\
    \hphantom{.....................................}
    mJ\Big[\big<JY,Z\big>JX-\big<JX,Z\big>JY\Big]= (mJ\hat S)(Z),
   \end{array}\\[2ex]
   \begin{array}{l}
   (\hat S S)(Z)=S_{JX\,JY}S_{XY}Z=S_{JX\,JY}\Big[\big<Y,Z\big>X
   -\big<X,Z\big>Y\Big]=\\
   \Big<JY,\big<Y,Z\big>X-\big<X,Z\big>Y\Big>JX-
    \Big<JX,\big<Y,Z\big>X-\big<X,Z\big>Y\Big>JY=\\
   \big<JY,X\big>\big<Y,Z\big>JX+\big<X,Z\big>\big<JX,Y\big>JY=\\
   \hphantom{........................................}
    mJ\Big[\big<Y,Z\big>X-\big<X,Z\big>Y\Big]=(mJ S)(Z),
   \end{array}\\[2ex]
   \begin{array}{l}
   (S J)(Z)=S_{XY}JZ=\big<Y,JZ\big>X-\big<X,JZ\big>Y=\\
   \hphantom{.....................................}
   J\big[\big<JY,Z\big>JX-\big<JX,Z\big>JY\big]=(J\hat S)(Z),
   \end{array}\\[2ex]
   \begin{array}{l}
   (\hat S J)(Z)=S_{JX\,JY}JZ=\big<JY,JZ\big>JX-\big<JX,JZ\big>JY=\\
   \hphantom{.....................................}
   J\big[\big<Y,Z\big>X-\big<X,Z\big>Y\big]=(J S)(Z),
   \end{array}
 \end{array}
    $$

    and the other entries of the table can be found in a similar way.

    From (\ref{tableC}) we see that $J(S+\hat S)=(S+\hat S)J$ and  we have
    $$
    \begin{array}{l}
    \ \ (S+\hat S)^2=S^2+\hat S^2+S\hat S+\hat SS=S^2+\hat S^2 +mJ(S+\hat S)
    \\[1ex]
    \begin{array}{ll}
    (S+\hat S)^3\,=&(S+\hat S)[S^2+\hat S^2 +mJ(S+\hat S)]=\\
    &S^3+\hat S^3+\hat SS^2+S\hat S^2+mJ(S+\hat S)^2=\\
    &-b^2S-b^2\hat S+(\hat SS)S+(S\hat S)\hat S+mJ(S+\hat S)^2=\\
    &-b^2(S+\hat S)+mJ(S^2 +\hat S^2)+mJ(S+\hat S)^2=\\
    &-b^2(S+\hat S)+mJ[(S+\hat S)^2-mJ(S+\hat S)]+\\
    & mJ(S+\hat S)^2=(m^2-b^2)(S+\hat S)+2mJ(S+\hat S)^2.
    \end{array}
    \end{array}
    $$
    Thus,
    \begin{equation}\label{qube1}
    (S+\hat S)^3=\Lin(S+\hat S, J(S+\hat S)^2)
    \end{equation}
    On the other hand, setting for brevity $R_{XY}=R$, from (\ref{curvop}) we derive
    \begin{equation}\label{qube2}
    \begin{array}{l}
    S+\hat S=\frac{4}{c}R-2mJ=\Lin(R,J),\\[1ex]
    (S+\hat S)^2=\Lin(R^2, JR, E).
    \end{array}
    \end{equation}
    Comparing (\ref{qube1}) and (\ref{qube2}) we conclude
    $$
    (S+\hat S)^3=\Lin\Big[\Lin(R,J),J\Lin(R^2,JR,E)\Big]=\Lin(JR^2,R,J).
    $$
    On the other hand, from $(\ref{qube2})_1$
    $$
    (S+\hat S)^3=\left(\frac{4}{c}\right)^3R^3+\Lin(JR^2,R,J).
    $$
    So, finally
    $$
    R^3=\Lin(JR^2,R,J).
    $$

    It is easy to trace that the coefficients of all linear combinations are
    polynomials in $\frac{1}{c}, b,m$. To complete the proof we should remark that
    $$
    \begin{array}{ll}
    R^4=R^3R=&\Lin(JR^2,R,J)R=\Lin(JR^3,R^2,JR)=\\
    &\Lin\Big[J\Lin(JR^2,R,J),R^2,JR\Big]=\Lin(R^2,JR,E)
    \end{array}
    $$
    which allows to find all powers of $R$ inductively.

    \pfend

    \textbf{Proof of Lemma \ref{Quat}}

    The curvature operator $R_{XY}$ of the quaternionic space form
   \\ $(M^{4n}(c);J_1,J_2,J_3;g)$ has the following expression
   $$
   \begin{array}{ll}
   R_{XY}Z=&\frac{\displaystyle c}{\displaystyle 4}\Big[\big<Y,Z\big>X-\big<X,Z\big>Y+\big<J_1Y,Z\big>J_1X-\big<J_1X,Z\big>J_1Y+\\
   &\big<J_2Y,Z\big>J_2X-\big<J_2X,Z\big>J_2Y+\big<J_3Y,Z\big>J_3X-\big<J_3X,Z\big>J_3Y+\\
   &2\,\big<X,J_1Y\big>J_1Z+2\,\big<X,J_2Y\big>J_2Z+2\,\big<X,J_3Y\big>J_3Z\Big].
   \end{array}
   $$
   where $J_1,\ J_2,\ J_3$ are operators of quaternionic structure
   $$
   J_1J_2=J_3, \, J_2J_3=J_1,\,
   \, J_3J_1=J_2, \,  J_i^2=-E,\,  \big<X,J_iY\big>=-\big<J_iX,Y\big>, \ i=\overline{1,3}
   $$
   Introduce the unit sphere type operator $S$ acting as
   $$
   S(Z)\stackrel{def}{=} S_{XY}Z=\big<Y,Z\big>X-\big<X,Z\big>Y,
   $$
   the operators $S_i(Z)$ acting as
   $$
   S_i(Z)\stackrel{def}{=} S_{J_iX\,J_iY}Z=\big<J_iY,Z\big>J_iX-\big<J_iX,Z\big>J_iY, \ i=\overline{1,3},
   $$
   and the operator $\hat S(Z)$ acting as
   $$
   \hat S(Z)\equiv S_1(Z)+S_2(Z)+S_3(Z).
   $$
   Finally, denote $m_i=\big<X,J_iY\big>,\ i=\overline{1,3}$, $m^2=m_1^2+m_2^2+m_3^2$,
    ${\cal J}=m_1J_1+m_2J_2+m_3J_3$. Then the curvature operator under
   consideration takes the form
   \begin{equation}\label{curvopq}
   R_{XY}Z=\frac{c}{4}\Big[ S+\hat S +2{\cal J}\Big]Z.
   \end{equation}
   Since $|X\wedge Y|=|(J_iX)\wedge (J_iY)| \ (i=\overline{1,3})$, the operators $S$ and $S_i$ satisfy
   $$
   S^3=-b^2S, \quad S_i^3=-b^2 S_i \ (i=\overline{1,3}),
   $$
   where $b^2=|X\wedge Y|^2$.

   The table of products for the operators $S$ and $\hat S$ is the following one.

   \begin{equation}\label{tableC2}
   \begin{array}{|c|c|c|c|c|c|c|c|} \hline
   & S & S_1 & S_2 & S_3 & J_1 & J_2 & J_3 \\[1ex] \hline
   S & S^2 & m_1J_1S_1 & m_2J_2S_2 & m_3J_3S_3 & J_1S_1 & J_2S_2 & J_3S_3 \\[1ex] \hline
   S_1 & m_1J_1S & S_1^2 & -m_3J_3S_2 & -m_2J_2S_3 & J_1S & J_2S_3 & J_3S_2 \\[1ex]   \hline
   S_2 & m_2J_2S & -m_3J_3S_1 & S_2^2 & -m_1J_1S_3 & J_1S_3 & J_2S & J_3S_1 \\[1ex] \hline
   S_3 & m_3J_3S & -m_2J_2S_1 & -m_1J_1S_2 & S_3^2 & J_1S_2 & J_2S_1 & J_3S \\[1ex] \hline
   J_1 & S_1J_1 & SJ_1 & S_3J_1 & S_2J_1 & -E & J_3 & -J_2 \\[1ex] \hline
   J_2 & S_2J_2 & S_3J_2 & SJ_2 & S_1J_2 & -J_3 & -E & J_1  \\[1ex] \hline
   J_3 & S_3J_3 & S_2J_3 & S_1J_3 & SJ_3 & J_2 & -J_1 & -E \\[1ex] \hline
\end{array}
   \end{equation}
   The expressions for products  $S S_i,\ S_i S,\ SJ_i$ one can find  similar to the table
   (\ref{tableC}) making formal replacements $\hat S\to S_i$ and $J\to J_i$. As
   concerns the other entries, we have
   $$
   \begin{array}{l}
   (S_1 S_2)(Z)=\\
   \qquad S_{J_1X\,J_1Y}
   S_{J_2X\,J_2Y}Z=S_{J_1X\,J_1Y}\Big[\big<J_2Y,Z\big>J_2X-\big<J_2X,Z\big>J_2Y\Big]=\\
   \qquad \Big<J_1Y,\big<J_2Y,Z\big>J_2X-\big<J_2X,Z\big>J_2Y\Big>J_1X-\\
   \qquad  \Big<J_1X,\big<J_2Y,Z\big>J_2X-\big<J_2X,Z\big>J_2Y\Big>J_1Y=\\
   \qquad
   \big<J_1Y,J_2X\big>\big<J_2Y,Z\big>J_1X+\big<J_1X,J_2Y\big>\big<J_2X,Z\big>J_1Y=\\
   \qquad
   J_1\Big[\big<X,J_3Y\big>\big<J_2Y,Z\big>X-\big<X,J_3Y\big>\big<J_2X,Z\big>Y\Big]=\\
   \qquad -J_1J_2\Big[m_3\big<J_2Y,Z\big>J_2X-m_3\big<J_2X,Z\big>J_2Y\Big]
   =(-m_3J_3S_2)(Z),\\
   \end{array}
  $$
  $$
   \begin{array}{l}
   (S_1 J_1)(Z)=S_{J_1X\,J_1Y}J_1Z=\big<J_1Y,J_1Z\big>J_1X-\big<J_1X,J_1Z\big>J_1Y=\\[1ex]
   \hphantom{.......................}J_1\big[\big<Y,Z\big>X-\big<X,Z\big>Y\big]=(J_1 S)(Z),
   \end{array}
   $$
   $$
   \begin{array}{l}
   (S_1 J_2)(Z)=S_{J_1X\,J_1Y}J_2Z=\big<J_1Y,J_2Z\big>J_1X-\big<J_1X,J_2Z\big>J_1Y=\\[1ex]
   \qquad J_1\big[\big<J_3Y,Z\big>X-\big<J_3X,Z\big>Y\big]=\\
   \hphantom{...........................} -J_1J_3\big[\big<J_3Y,Z\big>J_3X-\big<J_3X,Z\big>J_3Y\big]=
   (J_2 S_3)(Z)
   \end{array}
   $$
    and so on.

    From (\ref{tableC2}) we see that
    $$
    \begin{array}{l}
    (S+\hat S){\cal J}=(S+S_1+S_2+S_3)(m_1J_1+m_2J_2+m_3J_3)=\\
    m_1J_1S_1+m_2J_2S_2+m_3J_3S_3+m_1J_1S+
    m_2J_2S_3+m_3J_3S_2+m_1J_1S_3+\\
    m_2J_2S+m_3J_3S_1+m_1J_1S_2+m_2J_2S_1+m_3J_3S=\\
    (m_1J_1+m_2J_2+m_3J_3)(S+S_1+S_2+S_3)={\cal J} (S+\hat S).$$
    \end{array}
    $$

    Therefore, the operators $(S+\hat S)$ и ${\cal J}$
    commute and hence for the operator $\displaystyle R=\frac{c}{4}\{(S+\hat S)+2{\cal J}\}$
    the usual formula for powers can be applied:
    $$
    R^n=\left(\frac{c}{4}\right)^n\sum_{l=0}^n{n\choose l}\left(S+\hat S\right)^{n-l}
    2^l\left({\cal J} \right)^l.
    $$

    The powers for ${\cal J}$ can be found trivially, since
    $$
    \begin{array}{ll}
 {\cal J}^2=&m_1^2J_1^2+m_1m_2(J_1J_2+J_2J_1)+m_1m_3(J_1J_3+J_3J_1)+m_2^2J_2^2+\\[1ex]
 & m_2m_3(J_2J_3+J_3J_2)+m_3^2J_3^2=-m_1^2E-m_2^2E-m_3^2E=-m^2E,
    \end{array}
    $$
where $m_1^2+m_2^2+m_3^2=m^2$.

As concerns the powers of $S+\hat S$, the following proposition gives the answer.

\begin{proposition}The operator $S+\hat S$ possesses the recurrent property
$$
(S+\hat S)^5=-2\,(b^2+m^2)(S+\hat S)^3-(b^2-m^2)(S+\hat S),
$$
where $b^2=|X\wedge Y|^2$ and $
m^2=m_1^2+m_2^2+m_3^2=\big<X,J_1Y\big>^2+\big<X,J_2Y\big>^2+\big<X,J_3Y\big>^2$.
\end{proposition}
\begin{proof}
The proof is technical and in what follows we will use some auxiliary operator
products. Namely,

\begin{equation}\label{aux}
\begin{array}{l}
    \begin{array}{ll}
    S \hat S=S{\cal J}, & \hat S S={\cal J}S,\\
    S{\cal J}S=-m^2S, & S\hat S{\cal J}=-m^2S,\\
    S(S_1^2+S_2^2+S_3^2)=S^2{\cal J},& \hat S S^2={\cal J}S^2,\\
    \hat SS{\cal J}={\cal J}S{\cal J}, & \hat S{\cal J}S={\cal J}S^2,\\
    \end{array}
    \\[3ex]
    \begin{array}{l}
    \hat S^2=S_1^2+S_2^2+S_3^2-\hat S{\cal J}+{\cal J}S,\\
    \hat S(S_1^2+S_2^2+S_3^2)=-b^2\hat S-(S_1^2+S_2^2+S_3^2){\cal J}+{\cal
    J}S^2.
    \end{array}
\end{array}
  \end{equation}

    The proof is straightforward. Applying (\ref{tableC2}), we get

 \noindent
    $
    S \hat S=S(S_1+S_2+S_3)=m_1J_1S_1+m_2J_2S_2+
    m_3J_3S_3=\\
    \hphantom{.........} m_1SJ_1+m_2SJ_2+ m_3SJ_3=S{\cal J}.
    $

    In a similar way we find

    \noindent
    $
    \hat S S=(S_1+S_2+S_3)S=m_1J_1S+m_2J_2S+m_3J_3S={\cal J}S,
    $

    \noindent
    $
    \hat S^2=(S_1+S_2+S_3)(S_1+S_2+S_3)=S_1^2+S_2^2+S_3^2+S_1S_2+
    S_1S_3+S_2S_1+S_2S_3+S_3S_1+S_3S_2=S_1^2+S_2^2+S_3^2-
    m_3S_1J_3-m_2S_1J_2-m_3S_2J_3-m_1S_2J_2-m_2S_3J_2-
    m_1S_3J_1=S_1^1+S_2^2+S_3^2-\hat S{\cal J}+m_1S_1J_1+m_2S_2J_2+
    m_3S_3J_3=\\
    S_1^2+S_2^2+S_3^2-\hat S{\cal J} +{\cal J}S,
    $

\noindent\\
    $
    S{\cal J}S=S(m_1J_1+m_2J_2+m_3J_3)S=(m_1J_1S_1+m_2J_2S_2+
    m_3J_3S_3)S=-m_1^2S-m_2^2S-m_3^2S=-m^2S,
    $

\noindent\\
    $
    S\hat S{\cal J}= S{\cal J J}=-m^2S,
    $

\noindent\\
    $
    S(S_1^2+S_2^2+S_3^2)=m_1SJ_1S_1+m_2SJ_2S_2+m_3SJ_3S_3=m_1S^2J_1+
    m_2S^2J_2+m_3S^2J_3=S^2{\cal J},
    $

\noindent\\
    $
    \hat S{\cal J}S=(S_1+S_2+S_3)(m_1J_1+m_2J_2+m_3J_3)S=(m_1J_1S+
    m_2J_2S_3+m_3J_3S_2+m_1J_1S_3+m_2J_2S+m_3J_3S_1+m_1J_1S_2+
    m_2J_2S_1+m_3J_3S)S={\cal J}S^2+m_2J_2m_3J_3S+m_3J_3m_2J_2S+
    m_1J_1m_3J_3S+m_3J_3m_1J_1S+m_1J_1m_2J_2S+m_2J_2m_1J_1S={\cal
    J}S^2,
    $

\noindent\\
    $
    \hat S(S_1^2+S_2^2+S_3^2)=(S_1+S_2+S_3)(S_1^2+S_2^2+S_3^2)=
    S_1^3+S_2^3+S_3^3+S_1S_2^2+S_1S_3^2+S_2S_1^2+S_2S_3^2+S_3S_1^2+S_3S_2^2=
    -b^2\hat
    S-m_3J_3S_2^2-m_2J_2S_3^2-m_3J_3S_1^2-m_1J_1S_3^2-m_2J_2S_1^2-
    m_1J_1S_2^2=-b^2\hat S-m_3S_1^2J_3-m_2S_1^2J_2-m_3S_2^2J_3-m_1S_2^2J_1-
    m_2S_3^2J_2-m_1S_3^2J_1=-b^2\hat S-(S_1^2+S_2^2+S_3^2){\cal J}+m_1S_1^2J_1+
    m_2S_2^2J_2+m_3S_3^2J_3=-b^2\hat S+(S_1^2+S_2^2+S_3^2){\cal J}+m_1J_1S^2+
    m_2J_2S^2+m_3J_3S^2=\\
    -b^2\hat S +(S_1^2+S_2^2+S_3^2){\cal J}+{\cal J}S^2.
    $
\\

Now we are ready to find the powers of $(S+\hat S)$. Using (\ref{aux}), we get

    \noindent
    $
    (S+\hat S)^2=
    S^2+S\hat S+\hat SS+\hat S^2=S^2+S{\cal J}+{\cal J}S+S_1^2+S_2^2+S_3^2-
    \hat S{\cal J}+{\cal J}S=\\
    \hphantom{(.................}S^2+S{\cal J}+2{\cal J}S-\hat S{\cal J}+
    S_1^2+S_2^2+S_3^2.
    $

    Multiplying the result by $S+\hat S $ and applying again (\ref{aux}), we find

    \noindent
   $
    (S+\hat S)^3=
    (S+\hat S)[S^2+S{\cal J}+2{\cal J}S-\hat S{\cal J}+S_1^2+S_2^2+S_3^2]=
    S^3+S^2{\cal J}+2S{\cal J}S-S\hat S{\cal J}+S(S_1^2+S_2^2+S_3^2)+
    \hat SS^2+\hat SS{\cal J}+2\hat S{\cal J}S-\hat S^2{\cal J}+\hat S(S_1^2+S_2^2+S_3^2)=
    -b^2S+S^2{\cal J}-2m^2S+m^2S+S^2{\cal J}+{\cal J}S^2+{\cal J}S{\cal J}+
    2\hat S{\cal J}S-[(S_1^2+S_2^2+S_3^2)-\hat S{\cal J}+{\cal J}S]{\cal J}+[-b^2\hat S-
    (S_1^2+S_2^2+S_3^2){\cal J}+{\cal J}S^2]=-(b^2+m^2)S+2S^2{\cal J}+
    {\cal J}S^2+{\cal J}S{\cal J}+2{\cal J}S^2-(S_1^2+S_2^2+S_3^2){\cal J}+\hat S{\cal J}^2-
    {\cal J}S{\cal J}-b^2\hat S-(S_1^2+S_2^2+S_3^2){\cal J}+{\cal J}S^2=-(b^2+m^2)S+
    2S^2{\cal J}+4{\cal J}S^2-2(S_1^2+S_2^2+S_3^2){\cal J}-m^2\hat S-b^2\hat S=\\[1ex]
    \hphantom{(..............}
    -(b^2+m^2)(S+\hat S)+2S^2{\cal J}+4{\cal J}S^2-2(S_1^2+S_2^2+S_3^2){\cal J}.
    $

    Continue the process.

    \noindent
    $
    (S+\hat S)^4=
    (S+\hat S)[-(b^2+m^2)(S+\hat S)+2S^2{\cal J}+4{\cal J}S^2-
    2(S_1^2+S_2^2+S_3^2){\cal J}]=-(b^2+m^2)(S+\hat S)^2+2S^3{\cal J}+
    4S{\cal J}S^2-2S(S_1^2+S_2^2+S_3^2){\cal J}+2\hat SS^2{\cal J}+4\hat S{\cal J}S^2-
    2\hat S(S_1^2+S_2^2+S_3^2){\cal J}=-(b^2+m^2)(S+\hat S)^2-2b^2S{\cal J}-
    4m^2S^2-2S^2{\cal J}{\cal J}+2{\cal J}S^2{\cal J}+4{\cal J}S^3-2[-b^2\hat S-
    (S_1^2+S_2^2+S_3^2){\cal J}+{\cal J}S^2]{\cal J}=-(b^2+m^2)(S+\hat S)^2-
    2b^2S{\cal J}-4m^2S^2+2m^2S^2+2{\cal J}S^2{\cal J}-4b^2{\cal J}S+
    2b^2\hat S{\cal J}+2(S_1^2+S_2^2+S_3^2){\cal J}^2-2{\cal J}S^2{\cal J}=
    -(b^2+m^2)(S+\hat S)^2-2b^2S{\cal J}-2m^2S^2-4b^2{\cal J}S+
    2b^2\hat S{\cal J}-2m^2(S_1^2+S_2^2+S_3^2)=-(b^2+m^2)(S+\hat S)^2-
    2m^2[S^2+S{\cal J}+2{\cal J}S-\hat S{\cal J}+(S_1^2+S_2^2+S_3^2)]+
    (2m^2-2b^2)(S{\cal J}+2{\cal J}S-\hat S{\cal J})=\\[1ex]
    \hphantom{(..............}
    -(b^2+3m^2)(S+\hat S)^2+(2m^2-2b^2)(S{\cal J}+2{\cal J}S-\hat S{\cal
    J}).
    $

    Finally,

    \noindent
    $
    (S+\hat S)^5=(S+\hat S)[-(b^2+3m^2)(S+\hat S)^2+(2m^2-2b^2)(S{\cal J}+
    2{\cal J}S-\hat S{\cal J})]=-(b^2+3m^2)(S+\hat S)^3+
    (2m^2-2b^2)[S^2{\cal J}+2S{\cal J}S-S\hat S{\cal J}+\hat SS{\cal J}+2\hat S{\cal J}S-
    \hat S^2{\cal J}]=-(b^2+3m^2)(S+\hat S)^3+(2m^2-2b^2)[S^2{\cal J}-
    2m^2S-S{\cal J}^2+{\cal J}S{\cal J}+2{\cal J}S^2-(S_1^2+S_2^2+S_3^2-
    \hat S{\cal J}+{\cal J}S){\cal J}]=-(b^2+3m^2)(S+\hat S)^3+
    (2m^2-2b^2)[S^2{\cal J}-m^2S+{\cal J}S{\cal J}+2{\cal J}S^2-(S_1^2+S_2^2+
    S_3^2){\cal J}+\hat S{\cal J}^2-{\cal J}S{\cal J}]=-(b^2+3m^2)(S+\hat S)^3+
    (m^2-b^2)[2S^2{\cal J}+4{\cal J}S^2-2(S_1^2+S_2^2+S_3^2){\cal J}-2m^2S-
    2m^2\hat S]=-(b^2+3m^2)(S+\hat S)^3+(m^2-b^2)[2S^2{\cal J}+
    4{\cal J}S^2-2(S_1^2+S_2^2+S_3^2){\cal J}-(m^2+b^2)(S+\hat S)+
    (b^2-m^2)(S+\hat S)]=-(b^2+3m^2)(S+\hat S)^3+
    (m^2-b^2)[(S+\hat S)^3+(b^2-m^2)(S+\hat S)]=\\[1ex]
    \hphantom{(.................}
     -2(b^2+m^2)(S+\hat S)^3-(b^2-m^2)^2(S+\hat S)
    $

\noindent
    which completes the proof.
    \end{proof}
        Thus,
    \begin{equation}\label{Sfive}
    (S+\hat S)^5=\Lin((S+\hat S)^3, S+\hat S)
    \end{equation}
    On the other hand, setting for brevity $R_{XY}=R$, from (\ref{curvopq}) we derive
    \begin{equation}\label{qrec1}
     S+\hat S=\frac{4}{c}R-2{\cal J}=\Lin(R,{\cal J})
    \end{equation}
    Since $(S+\hat S)$ and ${\cal J}$ commute,  (\ref{curvopq}) implies
     the commutation of $R$ and ${\cal J}$. Keeping this and ${\cal J}^2=-m^2E$, from
     (\ref{qrec1})  we derive
     \begin{equation}\label{qrec3}
     (S+\hat S)^3=\left(\frac{4}{c}\right)^3R^3+\Lin({\cal J}R^2,R,{\cal J})
    \end{equation}
    \begin{equation}\label{qrec5}
     (S+\hat S)^5=\left(\frac{4}{c}\right)^5R^5+\Lin({\cal J}R^4,R^3,{\cal J}R^2,R,{\cal J})
    \end{equation}
    From (\ref{Sfive}), (\ref{qrec1}) and (\ref{qrec3})
    $$
    (S+\hat S)^5=\Lin\Big[\Lin(R^3,{\cal J}R^2,R,{\cal J}),\Lin(R,{\cal J})\Big]=
    \Lin(R^3,{\cal J}R^2,R,{\cal J}).
    $$
    So, finally from (\ref{qrec5})
    $$
    R^5=\Lin({\cal J}R^4,R^3,{\cal J}R^2,R,{\cal J}).
    $$

    It is easy to trace that the coefficients of all linear combinations are
    polynomials in $\frac{1}{c}, b,m$. To complete the proof we should remark that

    \noindent
    $
    R^6=R^5R=\Lin({\cal J}R^4,R^3,{\cal J}R^2,R,{\cal J})R=
    \Lin({\cal J}R^5,R^4,{\cal J}R^3,R^2,{\cal J}R)=\\
    \hphantom{.........}\Lin\Big[{\cal J}\Lin({\cal J}R^4,R^3,{\cal J}R^2,R,{\cal J}),R^4,{\cal J}R^3,R^2,{\cal
    J}R\Big]=\\
    \hphantom{.........}\Lin(R^4,{\cal J}R^3,R^2,{\cal J}R,E)
    $

    \noindent
    which allows to find all powers of $R$ inductively.

\pfend

\end{document}